\documentclass[suppldata]{interact}
\usepackage{graphicx}
\usepackage{setspace} 
\usepackage{multirow}
\usepackage{enumerate}
\usepackage{commath, amssymb,amsfonts,color,mathtools, enumerate}
\usepackage{hyperref}
\usepackage{amsmath}
\usepackage{anyfontsize}

\usepackage{epstopdf}
\usepackage[caption=false]{subfig}

\usepackage{cleveref} 
\usepackage{autonum} 

\usepackage{natbib}
\bibpunct[, ]{(}{)}{;}{a}{}{,}
\renewcommand\bibfont{\fontsize{10}{12}\selectfont}

\usepackage{url} 

\usepackage{afterpage}  

\begin{document}

\title{The skew-symmetric-Laplace-uniform distribution}

\author{
\name{Raju.~K. Lohot\textsuperscript{a} \href{https://orcid.org/0000-0003-0424-9447}{\includegraphics[scale=0.08]{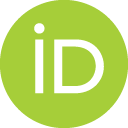}}  \thanks{CONTACT Raju.~K. Lohot. Email: rajulohot.92@gmail.com} and V. U. Dixit\textsuperscript{b} \href{https://orcid.org/0000-0002-4627-589X}{\includegraphics[scale=0.08]{orcid.png}}}
\affil{\textsuperscript{a}Department of Statistics, SVKM's Mithibai College of Arts, Chauhan Institute of Science \& Amrutben Jivanlal College of Commerce and Economics, Vile Parle (W), Mumbai, Maharashtra, India; \textsuperscript{b}Department of Statistics, University of Mumbai, Vidyanagari, Santacruz (E), Mumbai, Maharashtra, India}
}

\maketitle

\abstract{Laplace distribution is popular in the field of economics and finance. Still, data sets often show a lack of symmetry and a tendency of being bounded from either side of their support. In view of this, we introduce a new family of skew distribution using the skewing mechanism of \citet{azzalini1985class}, namely, skew-symmetric-Laplace-uniform distribution (SSLUD). Here uniform distribution is used not only to introduce skewness in Laplace distribution but also to restrict distribution support on one side of the real line. This paper provides a comprehensive description of the essential distributional properties of SSLUD. Estimators of the parameter are obtained using the method of moments and the method of maximum likelihood. The finite sample and asymptotic properties of these estimators are studied using simulation. It is observed that the maximum likelihood estimator is better than the moment estimator through a simulation study. Finally, an application of SSLUD to real-life data on the daily percentage change in the price of NIFTY 50, an Indian stock market index, is presented.}



\begin{amscode}
62E10, 62F10
\end{amscode}

\section{Introduction}

Symmetry is something which we try to seek naturally in everything, but not everything in the world is symmetric. So expecting symmetry in everything is unrealistic. In statistics, most classical procedures assume some kind of symmetry. However, the absence of symmetry is much more common in many data sets. In particular, much interest has been shown recently in a family of distributions called ``Skew-symmetric distributions". Let $f$ be a density function symmetric about zero, and $K$ an absolutely continuous distribution function such that the corresponding density function $K\, '$ is symmetric about zero. Then, Azzalini's form of skew-symmetric density function for any real $\lambda$, as mentioned in \citet{azzalini1985class}, is given as
\begin{equation} \label{azzalini}
    2\,f(x)\,K(\lambda x).
\end{equation}

\noindent \citet{arnold2004characterizations} studied a special case using $K$ as the cumulative distribution function (cdf) of $f$ in (\ref{azzalini}). \citet{nadarajah2003skewed} introduced the skew-symmetric-normal distribution family by replacing $f$ with $\phi$, the probability density function (pdf) of the standard normal distribution in (\ref{azzalini}). Further, they studied various skew-symmetric distributions by choosing $K$ as the cdf of normal, Student's t, Laplace, logistic, and uniform distributions. \citet{nadarajah2009skew} introduced and studied the skew logistic distribution considering $f$ and $K$ as pdf and cdf of logistic distribution, respectively in (\ref{azzalini}). 

When $f$ and $K$ are the density and distribution functions of the Laplace distribution  in (\ref{azzalini}), respectively, it is called a skew-Laplace distribution. \citet{aryal2005reliability} studied some properties of truncated skew-Laplace distribution, and \citet{kozubowski2008infinite} showed that a skew-Laplace distribution is infinitely divisible. Further, \citet{nekoukhou2012family} introduced a more general family of skew-Laplace distributions by considering $f$ as a standard Laplace pdf, $K$ as an arbitrary symmetric cdf, and $w$ as any odd continuous function in place of $\lambda x$ in (\ref{azzalini}). That is,
\begin{equation}
    e^{-\lvert x \rvert} \, F(w(x)).
\end{equation}

Recently much interest has been shown in the construction of flexible parametric classes of distributions that exhibit skewness and kurtosis, which is different from the normal distribution. While much of classical statistical analysis is based on Gaussian distributional assumptions, statistical modeling with the Laplace distribution has gained importance in many applied fields.  The motivation originates from data sets, including environmental, financial, and biomedical ones, which often do not follow the normal law. Models based on the Laplace distributions are popular in economics and finance; see \citet{zeckhauser1970linear, rachev1993laplace, ryden1998stylized, theodossiou1998financial, kotz2001laplace, kozubowski2001asymmetric}. They are rapidly becoming distributions of the first choice whenever “something” with heavier than normal tail is observed in the data. The interesting characteristic has often bound on its support from either end along with skew nature. i.e., data is positively skewed but bounded below or negatively skewed but bounded above. For example, consider the scenario of family income, which is typically positively skewed and bounded below by a certain amount. In this paper, by considering interesting applications of Laplace distribution, the need for skewness and restriction on the support of variable of interest, skew-symmetric-Laplace-uniform distribution (SSLUD) is introduced. Here, we consider $f$ as the standard Laplace density function and $K$ as a distribution function of Uniform$(-\theta, \theta)$ in (\ref{azzalini}). It provides a more flexible model representing the data as adequately as possible. Thus, we can expect this to be useful in more practical situations. The standard Laplace pdf is

\begin{equation} \label{lap}
    f(x)= \frac{1}{2} \,e^{-\lvert x \rvert}  , \quad x \in \mathbb{R}.
\end{equation}

\noindent The distribution function of Uniform$(-\theta, \theta)$ where $\theta >0$ is
\begin{equation}
K(x)  =
\begin{cases}
\; 0   & \text{if} \  x< -\theta, \\
\displaystyle\; \frac{x+\theta}{2\theta} & \text{if} \ - \theta \leqslant x < \theta, \\
\; 1 & \text{if} \ x \geqslant \theta.
\end{cases}
\end{equation}
   
\noindent Thus, the density function of SSLUD is 
\begin{equation} \label{gen.SSLLD}
    g(x)=2\, f(x) \, K(\lambda x) , \quad  x, \lambda \in \mathbb{R}.
\end{equation}

\noindent We define $\displaystyle \mu=\frac{\theta}{\lambda}$ so that model is identifiable. Here $\lambda \in \mathbb{R}$, $\theta > 0$ and hence $\mu \in \mathbb{R}-\{0\}$. Thus,
\begin{equation}\label{pdf.SSLUD}
    g(x) =
    \begin{cases}
    \; 0 & \text{if} \ \displaystyle \frac{x}{\mu} < -1 ,\\
    \displaystyle\; e^{-\lvert x \rvert} \ \left (\frac{x}{2\mu} + \frac{1}{2} \right)  & \text{if} \ \displaystyle -1 \leqslant \frac{x}{\mu} < 1,  \\
    \displaystyle\; e^{-\lvert x \rvert}\  & \text{if} \ \displaystyle \frac{x}{\mu} \geqslant 1.
    \end{cases} 
\end{equation}

\noindent  Here, one can notice that the support of $X$ is bounded above if $\mu < 0$ and bounded below if $\mu > 0$ by $-\mu$. The corresponding cdf $G(x)$ is as follows. When $\mu<0$,
\begin{subequations}
\begin{equation} \label{SSLUD.neg.cdf}
\begin{split}
G(x) & =
    \begin{cases}
      \displaystyle\; e^x & \text{if} \  x<\mu, \vspace{0.25 cm}\\
      \displaystyle\; \frac{e^{x}}{2\mu}(x+\mu-1)+\frac{e^{\mu}}{2\mu}&  \text{if} \ \mu \leqslant x<0,  \vspace{0.25 cm}\\
     \displaystyle\; 1+\frac{e^{\mu}}{2\mu}-\frac{e^{-x}}{2\mu}(x+\mu+1)&  \text{if} \ 0 \leqslant x< -\mu,  \vspace{0.25 cm}\\
    \; 1 &  \text{if} \  x \geqslant -\mu,
    \end{cases}\\
    \end{split}
    \end{equation}

\noindent \text{and when} $\mu>0$,
    
\begin{equation}
\label{SSLUD.postive.cdf}
\begin{split}
    G(x) & =
    \begin{cases}
       \; 0 & \text{if} \  x < -\mu, \vspace{0.25 cm} \\
      \displaystyle\; \frac{e^{x}}{2\mu}(x+\mu-1)+\frac{e^{-\mu}}{2\mu}&  \text{if} \ -\mu \leqslant x < 0,  \vspace{0.25 cm}\\
     \displaystyle\; 1+\frac{e^{-\mu}}{2\mu}-\frac{e^{-x}}{2\mu}(x+\mu+1)&  \text{if} \ 0 \leqslant x < \mu,  \vspace{0.25 cm}\\
    \; 1-e^{-x} &  \text{if} \  x \geqslant \mu.
    \end{cases}
\end{split}
\end{equation}
\end{subequations}

\begin{figure}
    \centering
    \includegraphics[width=0.8\textwidth]{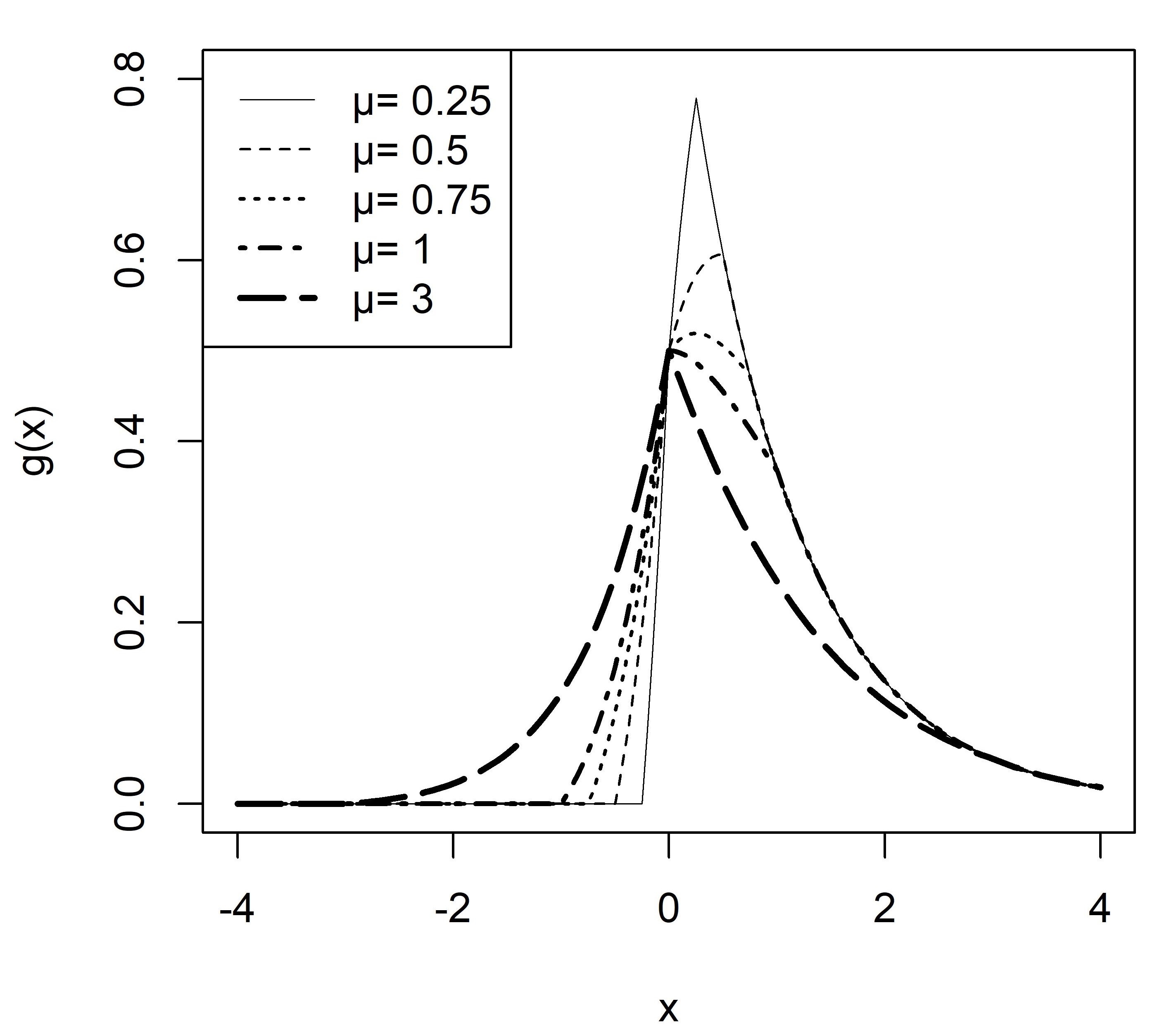}
    \caption{The skew-symmetric-Laplace-uniform pdf (\ref{pdf.SSLUD}) for $\mu=0.25, 0.5, 0.75, 1, 3$}
    \label{fig: pdf curve of SSLUD}
\end{figure}

\noindent Throughout the rest of this paper, unless otherwise stated, we shall assume that $\lambda > 0$, i.e., $\mu > 0$, since the corresponding results for $\lambda < 0$, i.e., $\mu < 0$, can be obtained using the fact that $-X$ has a pdf given by $2 f(x) K(-\lambda x)$. Figure \ref{fig: pdf curve of SSLUD} illustrates the shape of the pdf (\ref{pdf.SSLUD}) for $\mu=0.25, 0.5, 0.75, 1, 3$. 

The skew-symmetric-Laplace-uniform distribution with parameter $\mu$, $SSLUD(\mu)$ appears not to have been introduced yet. We provide a comprehensive description of the mathematical properties of (\ref{pdf.SSLUD}). This paper follows up on \citet{nadarajah2009skew}, where a comprehensive description of the mathematical properties for the skew-logistic distribution is provided. Here, we have derived formulas for moment generating function, characteristic function, and first four raw moments (Sect. \ref{Section MGF, CGF, and moments}), mode and median (Sect. \ref{section Mode and Median}), hazard rate function (Sect. \ref{sect hazard rate function}), mean deviation about `$a$' (Sect. \ref{sect Mean Deviation}), R\`enyi entropy and Shannon entropy (Sect. \ref{sect entropy}), simulation and estimation by the methods of moments and maximum likelihood (Sect. \ref{sect estimation}). We also discuss these estimators’ finite sample and asymptotic properties (Sect. \ref{sect estimation}). Finally, the application of $SSLUD(\mu)$ to real-life data on the daily percentage change in the price of NIFTY 50, an Indian stock market index, is discussed. Comparison of fitting of $SSLUD(\mu)$ is done with fitting of normal distribution $N(\theta, \sigma^2)$, Laplace distribution $L(\theta, \beta)$, and skew-Laplace distribution $SL(\lambda)$ for the above data (Sect. \ref{sect application}).

\section{Moment generating function, characteristic function, and moments} \label{Section MGF, CGF, and moments}
Here, we derive the moment generating function and the characteristic function of r. v. $X$ having pdf given in (\ref{pdf.SSLUD}). The moment generating function (MGF) is $M_{X}(t) = E(e^{tX})$. By using (\ref{pdf.SSLUD}), one obtains

\begin{equation}
      M_X(t)=\frac{1}{2 \mu} \left \{ \frac{-1+e^{- \mu (1+t)}}{(1+t)^2} + \frac{1-e^{-\mu (1-t)}}{(1-t)^2}\right \} + \frac{1}{(1-t^2)}\ , \quad \text{for} \ t<1.
\end{equation}

\noindent The corresponding characteristic function defined by $\phi_X(t) = E(e^{itX})$ is given as 

\begin{equation}
      \phi_X(t)=\frac{1}{2 \mu} \left \{ \frac{-1+e^{- \mu (1+it)}}{(1+it)^2} + \frac{1-e^{-\mu (1-it)}}{(1-it)^2}\right \} + \frac{1}{(1+t^2)} \ , \quad \text{for} \ it<1,
\end{equation}

 \noindent where $i = \sqrt{-1}$ is the complex imaginary unit.
 
The moments of a probability distribution are a collection of descriptive constants used for measuring its properties. Here, we derive the expression of the first four raw moments of $X$. They are as follows.
\begin{equation} \label{raw}
    \begin{split}
        \mu_{1}^{'} & = \frac{2}{\mu}-\left(1+\frac{2}{\mu}\right)e^{-\mu}\, ,\\
        \mu_{2}^{'} & =2 \, ,\\
        \mu_{3}^{'} & =\frac{24}{\mu}- e^{-\mu} \left[\mu^2+6\mu+18+\frac{24}{\mu}\right] \, ,\\
        \mu_{4}^{'} & =24 \, .
    \end{split}
\end{equation}

\noindent We see that $\mu_{2r}^{'}=(2r)!$ for $r=1, 2, \ldots$ and corresponding central moments can be obtained using these raw moments but can not be simplified further. Note that, expressions given in (\ref{raw}) are valid only for $\mu>0$. If $\mu<0$, one must replace $\mu$ by $-\mu$ in each of these expressions; in addition, the expressions for the odd order moments must be multiplied by -1.

\begin{figure}
    \centering
    \includegraphics[width=0.8\textwidth]{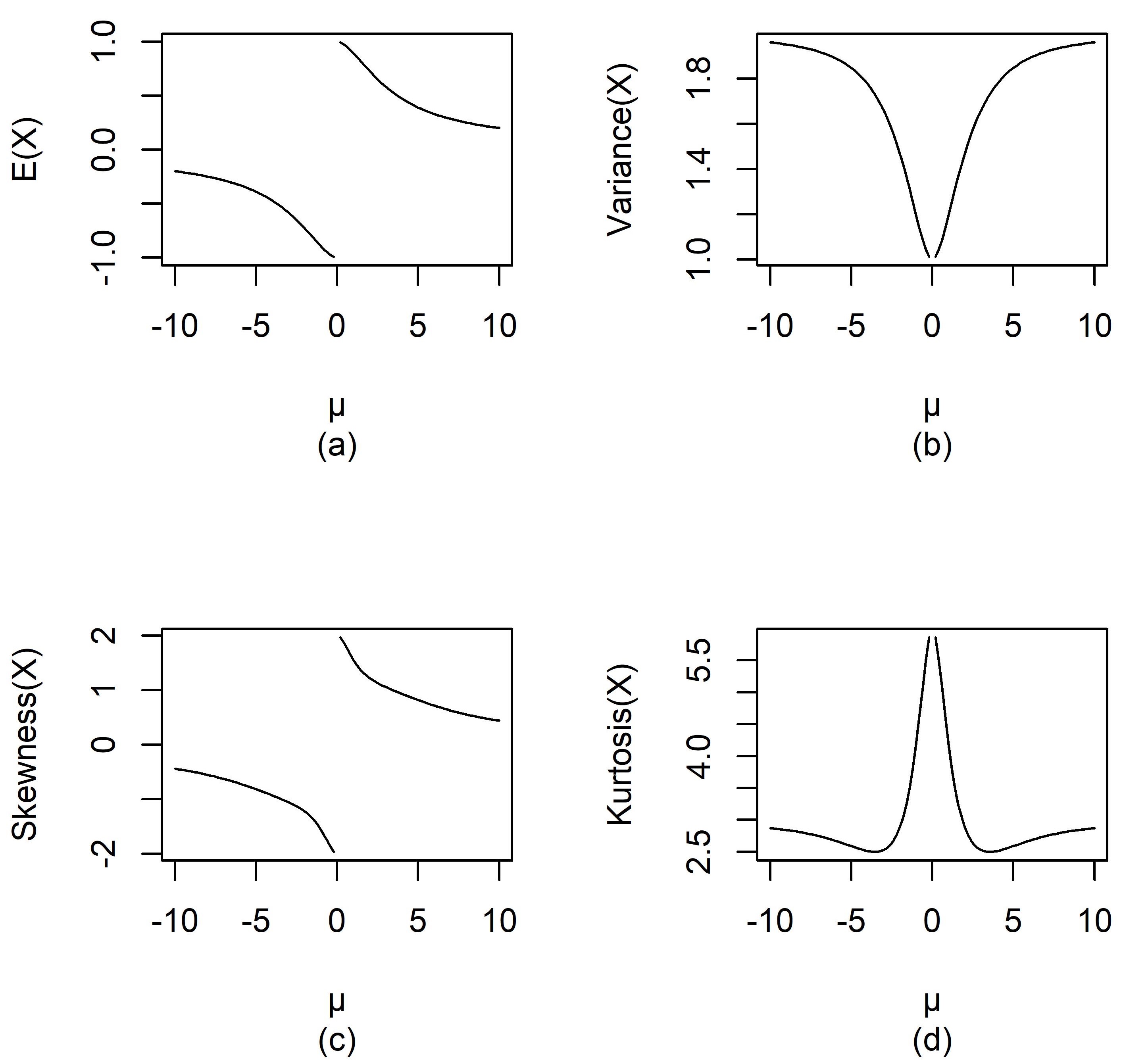}
    \caption{Variation of the four measures: (a) E($X$), (b) Variance($X$), (c) Skewness($X$) and \\ (d) Kurtosis($X$) for $\mu = -10, \ldots, 10$}
    \label{fig: plot of mean, var, skew and kur}
\end{figure}

Figure \ref{fig: plot of mean, var, skew and kur} illustrates the behavior of the four measures E($X$), Var($X$), Skewness($X$) and Kurtosis($X$) for $\mu=-10,\ldots,10$. Mean and skewness are decreasing functions of $\mu$ over the range $(-\infty,0)$ and $(0, \infty)$, while variance and kurtosis are even functions of $\mu$. The variance strictly decreases as $\mu$ moves from $-\infty$ to 0 and increases as $\mu$ moves from 0 to $\infty$.

\section{Mode and median} \label{section Mode and Median}

Mode is the value of the r. v. $X$ at which pdf $g(x)$ is maximum. When $\mu >0$, $g'(x)$ is,

\begin{equation} \label{g dash x}
    g'(x)= 
    \begin{cases} 
    \displaystyle\; e^{x} \left( \frac{x+1}{2\mu}+ \frac{1}{2}\right) & \text{if} \  -\mu \leqslant x < 0,\\
    \displaystyle\; e^{-x} \ \left (\frac{1-x}{2\mu} - \frac{1}{2} \right)  & \text{if} \ 0 \leqslant x < \mu,  \\
    \displaystyle\; - \, e^{-x}\  & \text{if} \  x \geqslant \mu.
    \end{cases}
\end{equation}

\noindent It is clear from (\ref{g dash x}) that the function $g(x)$ is increasing in $[-\mu, 0)$ and decreasing in $[\mu, \infty)$. Hence, mode $M_0$ of (\ref{pdf.SSLUD}) lies in the interval $[0, \mu]$. Accordingly, we equate $g'(x)$ to zero and solve for x. Thus, the value of $M_0$ is $M_0 = 1-\mu$ for $\displaystyle \frac{1}{2} \leqslant \mu < 1$. But when $\mu < \displaystyle \frac{1}{2}$, the function $g(x)$ increases in $[-\mu, \mu)$ and decreases in $[\mu, \infty)$. Hence, $M_0 = \mu$. Similarly, when $\mu > 1$, the function $g(x)$ increases in $[-\mu, 0)$ and decreases in $[0, \infty)$. Hence, $M_0 = 0$. On similar lines, one can derive the expression of $M_0$ for $\mu<0$. Thus, combining these two expressions of $M_0$, we get 

\begin{equation}
    M_0= 
    \begin{cases} 
    \; \mu & \text{if}  \ \displaystyle 0 < |\mu| < \frac{1}{2},\\
    \; \text{sign}(\mu)-\mu & \text{if} \ \displaystyle \frac{1}{2} \leqslant |\mu| < 1,\\
    \; 0 &  \text{if} \ |\mu| \geqslant 1.
    \end{cases}
\end{equation}

The median $M$ of (\ref{pdf.SSLUD}) is the value of r. v. $X$ such that $G(M)=\frac{1}{2}$. Thus, for $\mu > 0$ using (\ref{SSLUD.postive.cdf}),  

\begin{equation}
   M =
    \begin{cases}
      \; \text{Solution of the equation,} & \text{if} \  \displaystyle G(0)>\frac{1}{2}, \\  \; e^{M}(M-1+\mu)+e^{-\mu}-\mu=0   \vspace{0.35 cm} \\ 
      \; \text{Solution of the equation,} & \text{if} \ \displaystyle G(0) \leqslant \frac{1}{2} < G(\mu),\\ \; e^{-M}(-M-1-\mu)+e^{-\mu}+\mu=0  \vspace{0.35 cm} \\ 
      \; \ln 2  & \text{if} \ \displaystyle G(\mu) \leqslant \frac{1}{2},  
      
    \end{cases}         
\end{equation}

\noindent where  $\displaystyle G(0)=\frac{1}{2}+\frac{e^{-\mu}-1}{2\mu}$ and $\displaystyle G(\mu)=1-e^{-\mu}$. Table \ref{Table of median} represents values of the median of (\ref{pdf.SSLUD}) for different positive values of $\mu$ using the Newton-Raphson iterative procedure in R-software. If $\mu<0$, one can obtain the median $M$ on similar lines using (\ref{SSLUD.neg.cdf}).
    
\begin{table}
\caption{Median of SSLUD($\mu$) for $\mu=$0.25, 0.5, \ldots , 1.5}\label{Table of median}%
\begin{tabular}{@{}lcccccc@{}}
\toprule
$\mu$ & 0.25 & 0.5 & 0.75 & 1 & 1.25 & 1.5 \\
\midrule
$M$ & 0.6931472 & 0.6931472 & 0.6920484 & 0.6681079 & 0.6273646 & 0.5811654 \\
\bottomrule
\end{tabular}
\end{table}

\section{Hazard rate function} \label{sect hazard rate function}
The reliability function $R(x)=1-G(x)$ for $\mu > 0$ is obtained using (\ref{SSLUD.postive.cdf}) as,
\begin{equation}
R(x) =
    \begin{cases}
     \; 1 & \text{if} \  x < -\mu, \vspace{0.25 cm} \\
     \displaystyle  \; 1-\frac{e^{x}}{2\mu}(x+\mu-1)-\frac{e^{-\mu}}{2\mu}&  \text{if} \ -\mu \leqslant x < 0,  \vspace{0.25 cm}\\
     \displaystyle\; -\frac{e^{-\mu}}{2\mu}+\frac{e^{-x}}{2\mu}(x+\mu+1)&  \text{if} \ 0 \leqslant x < \mu,  \vspace{0.25 cm}\\
   \displaystyle\; e^{-x} &  \text{if} \  x \geqslant \mu.
    \end{cases}
\end{equation}

\noindent The hazard rate function is an important quantity,  characterizing life phenomena. After some simple steps, one can get the hazard function $\displaystyle h(x)=\frac{g(x)}{R(x)}$ for $\mu>0$ as follows.

\begin{equation}
    h(x) =
    \begin{cases}
      \; 0 & \text{if} \  x < -\mu, \vspace{0.25 cm} \\
      \displaystyle\; \left[-1+\frac{1+(2\mu-e^{-\mu}) e^{-x}}{x+\mu}\right]^{-1}&  \text{if} \ -\mu \leqslant x < 0,  \vspace{0.25 cm}\\
      \displaystyle\; \left[1+\frac{1-e^{(x-\mu)}}{x+\mu}\right]^{-1}&  \text{if} \ 0 \leqslant x < \mu,  \vspace{0.25 cm}\\
    \; 1 &  \text{if} \  x \geqslant \mu.
    \end{cases}
\end{equation}

\noindent One can easily check that $h(x)$ is increasing function of $x$ for $\mu<0$ as well as for $\mu>0$. Hence, $SSLUD(\mu)$ is increasing failure rate (IFR) distribution.

\section{Mean deviation} \label{sect Mean Deviation}
The amount of scatter in a population is evidently measured to some extent by the totality of deviations from the mean and median. These are known as the mean deviation about the mean and the mean deviation about the median, respectively. Mean deviation about an arbitrary real number `$a$' is defined by $\eta_a=E\lvert X-a \rvert$.

\noindent It leads to expression as 
\begin{equation}
    \eta_a =
    \begin{cases}
    \displaystyle\; -\left(1+\frac{2}{\mu}\right) e^{-\mu}+\left(\frac{2}{\mu}-a\right) & \text{if} \  a<-\mu, \vspace{0.35cm}\\    
    \displaystyle\; \left(\frac{a}{\mu}\right)e^{-\mu}+\left(\frac{a}{\mu}-\frac{2}{\mu}+1\right)e^{a}+\left(\frac{2}{\mu}-a\right) & \text{if} \  -\mu\leqslant a < 0, \vspace{0.35 cm}\\
   \displaystyle\;  \left(\frac{a}{\mu}\right)e^{-\mu}+\left(\frac{a}{\mu}+\frac{2}{\mu}+1\right)e^{-a}+\left(a-\frac{2}{\mu}\right) & \text{if} \  0\leqslant a < \mu, \vspace{0.35cm}\\
   \displaystyle\;  \left(1+\frac{2}{\mu}\right)e^{-\mu}+2e^{-a}+\left(a-\frac{2}{\mu}\right) & \text{if} \  a \geqslant \mu.
    \end{cases}
\end{equation}

\noindent To obtain mean deviation about mean and mean deviation about median, `$a$' in the above expression can be replaced by mean and median, respectively.

\section{Entropy} \label{sect entropy}
The entropy of a random variable $X$ measures the variation of uncertainty. The R\`enyi entropy of order $\alpha$ is 
\begin{equation}\label{renyi.def}
    H_{\alpha}=\frac{1}{1-\alpha} \, \log_{2} \left\{\int g^{\alpha}(x)dx\right\}, \ \ \alpha > 0 \  , \ \alpha \neq 1,
\end{equation}

\noindent where $g(x)$ is pdf of random variable $X$. By using (\ref{pdf.SSLUD}), one can write 
    \[ \int g^{\alpha}(x)dx = I_1 + I_2 + \frac{e^{-\alpha x}}{\alpha}, \]
\noindent where \( I_1=\int\displaylimits_{-\mu}^{0} e^{\alpha x} \left(\frac{x}{2\mu}+\frac{1}{2}\right)^\alpha dx 
    =\left(\frac{1}{\alpha}\right)\left(\frac{1}{2}\right)^\alpha \left[\sum_{j=0}^{\alpha}\left(\frac{-1}{\mu \alpha}\right)^j \frac{\alpha!}{(\alpha-j)!}-\left(\frac{-1}{\mu \alpha}\right)^\alpha e^{-\mu \alpha}\right] \)
and \( I_2=\int\displaylimits_{0}^{\mu} e^{-\alpha x} \left(\frac{x}{2\mu}+\frac{1}{2}\right)^\alpha dx 
    =\sum_{j=0}^{\alpha} \frac{\alpha!}{(\alpha-j)! (2\mu \alpha)^j} \left[\frac{2^{-(\alpha-j)}-e^{-\mu \alpha}}{\alpha}\right]. \) 
    
Therefore,
\begin{equation} \label{renyi.int.fx}
\begin{split}
     \int g^{\alpha}(x)dx= &  \frac{1}{\alpha} \sum_{j=0}^{\alpha} \frac{\alpha!}{(\alpha-j)!} \left[\frac{2^{-(\alpha-j)}(1+(-1)^j)-e^{-\mu \alpha}}{(2\mu\alpha)^j}\right]\\
     & + \frac{e^{-\mu\alpha}}{\alpha}\left[1-\left(\frac{-1}{2\mu\alpha}\right)^\alpha \right].
\end{split}
\end{equation}

\noindent One can obtain the R\`enyi entropy of order $\alpha$ by substituting (\ref{renyi.int.fx}) in (\ref{renyi.def}).

The Shannon entropy function is the particular case of (\ref{renyi.def}) for $\alpha \uparrow 1$, and it is $H=E[-\log_{2} g(X)]$, where $g(x)$ is pdf of random variable $X$. Using this definition, after some simplification we get,
\begin{equation}
        \begin{split}
            H= & \frac{1}{\ln 2}-\int\displaylimits_{0}^{\mu} \frac{xe^{-x}}{2\mu} \  \log_2 \left[\frac{\left(\displaystyle \frac{x}{2\mu}+\frac{1}{2}\right)}{\left(\displaystyle \frac{-x}{2\mu}+\frac{1}{2}\right)}\right] dx\\
            & -\int\displaylimits_{0}^{\mu} \frac{e^{-x}}{2} \  \log_2 \left[\left(\frac{-x}{2\mu}+\frac{1}{2}\right)\left(\frac{x}{2\mu}+\frac{1}{2}\right)\right] dx.
        \end{split}
\end{equation}

\begin{figure}
    \centering
    \includegraphics[width=0.8\textwidth]{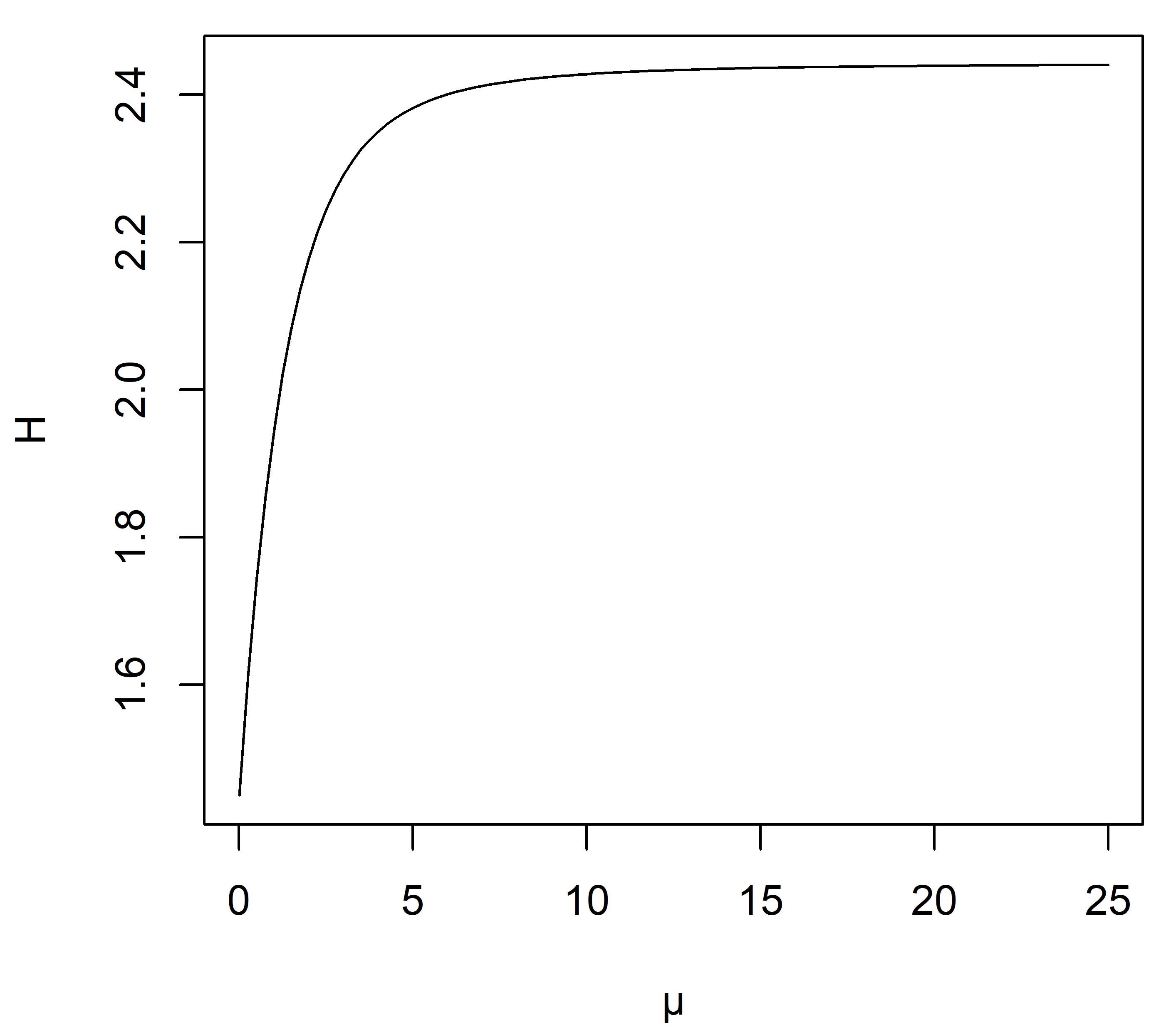}
    \caption{Behavior of Shannon entropy function}
    \label{fig: Behavior of Shannon entropy function}
\end{figure}

\noindent Since the above integration is cumbersome, we numerically evaluate $H$ for different values of $\mu$ using R-software. Figure \ref{fig: Behavior of Shannon entropy function} represents a graph of $\mu$ $(\mu > 0)$ versus $H$.

\section{Estimation} \label{sect estimation}
Here, we first consider simulating values of a random variable $X$ with the pdf (\ref{pdf.SSLUD}) using the inverse transformation technique. Let $r$ be a random number between zero and one. 
The generator to generate a random sample is 

\begin{equation} \label{generator}
 X =
    \begin{cases}
      \; \text{Solution of the equation,} & \text{if} \  0\leqslant r < G(0),\\  \; e^{x}(x-1+\mu)+e^{-\mu}-2r\mu=0   \vspace{0.35 cm} \\ 
      \; \text{Solution of the equation,} & \text{if} \  G(0) \leqslant r < G(\mu),\\ \; e^{-x}(-x-1-\mu)+e^{-\mu}+2(1-r)\mu=0  \vspace{0.35 cm} \\ 
      \; -\ln (1-r)  & \text{if} \  G(\mu) \leqslant r \leqslant 1.
      
    \end{cases}       
\end{equation}

\noindent One can use the Newton-Raphson method to solve the equation in (\ref{generator}) and generate a random sample from $SSLUD(\mu)$ given in (\ref{pdf.SSLUD}).  

Now, we consider the estimation of $\mu$ by the method of moments and the method of maximum likelihood. To estimate unknown parameter $\mu$, we have to consider both the cases $\mu<0$ and $\mu>0$ together. Suppose $x_1,...,x_n$ is an observed random sample of size `$n$' from (\ref{pdf.SSLUD}). For the method of moments estimation, after equating sample mean $\overline{x}$ to the first population raw moment of (\ref{pdf.SSLUD}), one obtains the equation 
\begin{equation} \label{mom.eq.1}
   \overline{x} =
   \begin{cases}
   \displaystyle\; \frac{2}{\mu}+\left(1-\frac{2}{\mu}\right)e^{\mu} & \quad \text{if} \  \mu < 0,\\
   \displaystyle\; \frac{2}{\mu}-\left(1+\frac{2}{\mu}\right)e^{-\mu} & \quad \text{if} \ \mu > 0.
   \end{cases}
\end{equation}

\noindent From Figure \ref{fig: plot of mean, var, skew and kur}, we see that $\mu_1^{'}$ decreases from 0 to -1 when $-\infty<\mu<0$ and it decreases from 1 to 0 when $0<\mu<\infty$, i.e., always $-1<\mu_1^{'}<1$. Therefore, if $\overline{x}<-1$ or $\overline{x}>1$ for a particular sample, then (\ref{mom.eq.1}) will not have an exact solution. As per Figure \ref{fig: plot of mean, var, skew and kur}, $\mu$ corresponds to the closest value of $\mu_1^{'}$ to $\overline{x}$ if $\overline{x}<-1$ or $\overline{x}>1$ is a value close to zero. But as per parameter space, $\mu$ can not take the value zero. Hence, we define the moment estimator $\tilde{\mu}$ of $\mu$ as $-10^{-5}$ if $\overline{x}<-1$ and $10^{-5}$ if $\overline{x}>1$. Thus, the moment estimator $\tilde{\mu}$ of $\mu$ is obtained as follows.

\begin{equation}
        \tilde{\mu}  = 
        \begin{cases}
        \; -10^{-5} & \quad \text{if} \ \  \overline{x} \leqslant -1,\vspace{0.35 cm} \\
        \; \text{Solution of the equation,}  & \quad \text{if} \  \ -1<\overline{x} < 0, \\
        \; \displaystyle \frac{2}{\mu}+\left(1-\frac{2}{\mu}\right)e^{\mu}-\overline{x} = 0 \vspace{0.35 cm} \\
        \; \text{Solution of the equation,}  & \quad \text{if} \ \ 0\leqslant \overline{x}<1, \\
        \; \displaystyle \frac{2}{\mu}-\left(1+\frac{2}{\mu}\right)e^{-\mu}-\overline{x}=0 \vspace{0.35 cm} \\
        \; 10^{-5} & \quad \text{if} \ \ \overline{x} \geqslant 1.
        \end{cases}
\end{equation}

We consider the estimation of $\mu$ by the method of maximum likelihood in the following. Let $x_{(1)}, x_{(2)},\ldots, x_{(n)}$ be the order statistics of given sample. Suppose $\mu<0$ and `$r_1$' denotes the number of observations less than $\mu$ such that $-\infty<x_{(1)}<x_{(2)}<\ldots < x_{(r_1)}\leqslant\mu \leqslant x_{(r_1 +1)}<\ldots<x_{(n)}<-\mu<\infty$\textcolor{red}{,} i.e. $-\infty<\mu< \min(0,\, -x_{(n)})$ where $r_1=0,1,2,\ldots, n$. Similarly, suppose $\mu>0$ and `$r_2$' denotes the number of observations lying in the interval $[-\mu, \mu]$ such that $-\infty<-\mu<x_{(1)}<x_{(2)}<\ldots < x_{(r_2)} \leqslant \mu \leqslant x_{(r_2 +1)} <\ldots<x_{(n)}< \infty$\textcolor{red}{,} i.e. $\max(0,\, -x_{(1)}) < \mu < \infty$ where $r_2=0,1,2,\ldots, n$. Hence, the log-likelihood function of $\mu$ is written as
\begin{equation}
    l=
    \begin{cases}
    \; l_1 = \displaystyle -\sum\limits_{i=1}^{n} \lvert x_{(i)}\rvert+\sum\limits_{i=r_1 +1}^{n} ln\left(\frac{x_{(i)}+\mu}{2\mu}\right) & \text{if} \ -\infty<\mu < \min\{0, \, -x_{(n)}\},\\
    \; l_2 =\displaystyle -\sum\limits_{i=1}^{n} \lvert x_{(i)} \rvert+\sum\limits_{i=1}^{r_2} ln\left(\frac{x_{(i)}+\mu}{2\mu}\right) & \text{if} \ \max\{-x_{(1)}, \, 0\}<\mu<\infty.
    \end{cases}
\end{equation}

\noindent In the following, we give a step-wise procedure for computation of the MLE $\hat{\mu}$ of $\mu$.
 \begin{enumerate}[\bfseries Step 1:]
     \item Numerically maximize $l_1$ over the range $(-a, \min\{0, -x_{(n)}\})$. Suppose the maximum value of $l_1$ is $\hat{l}_1$ which is attained at $\hat{\mu}_1$, say, where `$a$' is a sufficiently large positive number chosen for computation purposes.
     
     \item Numerically maximize $l_2$ over the range $(\max\{-x_{(1)},0\}, a)$. Suppose the maximum value of $l_2$ is $\hat{l}_2$ which is attained at $\hat{\mu}_2$, say.

     \item MLE $\hat{\mu}$ of $\mu$ is \begin{equation}
    \begin{aligned}
    \hat{\mu} = & 
    \begin{cases}
    \; \hat{\mu}_1 \quad \text{if} \ \hat{l}_1 > \hat{l}_2,\\
    \; \hat{\mu}_2  \quad \text{otherwise}.
    \end{cases}
    \end{aligned}
    \end{equation}
 \end{enumerate}

Finite sample properties of $\tilde{\mu}$ and $\hat{\mu}$ are studied using simulation, and computations are done using R- software. Table \ref{Table of Bias and MSE 1} and Table \ref{Table of Bias and MSE 2} presents bias and MSE of $\tilde{\mu}$ and $\hat{\mu}$ for $n=100 (100) 1000$ and for $\mu = -1.5, -0.75, -0.25, 0.25, 0.75, 1.5$. We see that bias and MSE decrease as sample size $n$ increases for both MLE $\hat{\mu}$ and moment estimator $\tilde{\mu}$, with few exceptions only for bias. Further, the MSE of $\hat{\mu}$ is always less than the corresponding MSE of $\tilde{\mu}$. Also, one can observe that sign of bias of MLE $\hat{\mu}$ is opposite to the sign of parameter $\mu$. As parameter $\mu$ approaches zero from any side, MSE and magnitude of bias of $\hat{\mu}$ decrease. But, no such observation in the case of the moment estimator $\tilde{\mu}$. To check the asymptotic nature of the distribution of $\hat{\mu}$ and $\tilde{\mu}$ using simulation, we plotted observed densities for various values of the sample size $n$. We observe that as $n$ increases, the distribution of both $\hat{\mu}$ and $\tilde{\mu}$ converges to the normal distribution, but the rate of convergence to normal distribution seems to be much higher for $\hat{\mu}$ than $\tilde{\mu}$. Thus, based on all the above results, we conclude that MLE is better than the moment estimator of $\mu$ for SSLUD($\mu$).

\begin{table}
\caption{Bias and MSE of MLE and moment estimator for $\mu= -1.5, -0.75, -0.25$, sample size $n=100(100)1000$, and simulation size $N=2000$}
\renewcommand{\arraystretch}{1}
\begin{tabular*}{\textwidth}{@{\extracolsep\fill}cccccc} 
\toprule
\multirow{2}{1 cm}{$\mu$} & \multirow{2}*{$n$} & \multicolumn{2}{@{}c@{}}{MLE}  & \multicolumn{2}{@{}c@{}}{Moment estimator}\\ \cmidrule{3-4}\cmidrule{5-6}%
& & Bias & MSE & Bias & MSE \\ \midrule

\multirow{10}{1 cm}{-1.5}& 100 & 0.06024381 & 0.045176654 & 0.0003797204 & 0.50541564\\
& 200 & 0.03906522 & 0.021357856 & 0.0157055747 & 0.26786436\\
& 300 & 0.02933963 & 0.012989993 & -0.0008644177 & 0.16465504\\
& 400 & 0.02375948 & 0.009209691 & -0.0085926100 & 0.12896987\\
& 500 & 0.01766114 & 0.007451091 & -0.0065899772 & 0.10033671\\
& 600 & 0.01808905 & 0.006370478 & 0.0075580998 & 0.08311324\\
& 700 & 0.01456689 & 0.005214147 & -0.0003676365 & 0.06853064\\
& 800 & 0.01492506 & 0.004456640 & -0.0164346584 & 0.05661172\\
& 900 & 0.01381575 & 0.003861085 & -0.0079539753 & 0.05188121\\
& 1000 & 0.01378036 & 0.003294051 & 0.0052817371 & 0.04687188\\ 
\midrule
\multirow{10}{1 cm}{-0.75} & 100 & 0.032511901 & 0.0134889603 & -0.02488011 & 0.37927153\\
& 200 & 0.023375017 & 0.0062341973 & 0.02678563 & 0.24807994\\
& 300 & 0.022614234 & 0.0040753483 & 0.05442608 & 0.18159122\\
& 400 & 0.013782428 & 0.0026481775 & 0.02407104 & 0.14540491\\
& 500 & 0.014188578 & 0.0019942288 & 0.04011415 & 0.12661616\\
& 600 & 0.010531113 & 0.0016580024 & 0.02741234 & 0.10965835\\
& 700 & 0.009567891 & 0.0014142729 & 0.02247448 & 0.09215384\\
& 800 & 0.008969822 & 0.0012107034 & 0.03299546 & 0.08416938\\
& 900 & 0.008833928 & 0.0010258905 & 0.02604549 & 0.07738691\\
& 1000 & 0.007307065 & 0.0009473139 & 0.01597401 & 0.06532169\\ 
\midrule
\multirow{10}{1 cm}{-0.25} & 100 & 0.026910685 & 0.0043941552 & -0.20148295 & 0.28701317\\
& 200 & 0.014633591 & 0.0016442189 & -0.12561277 & 0.18340413\\
& 300 & 0.011298482 & 0.0010134272 & -0.11008763 & 0.14648063\\
& 400 & 0.008403715 & 0.0007178242 & -0.07115405 & 0.11693260\\
& 500 & 0.006866262 & 0.0005193217 & -0.05475621 & 0.10104767\\
& 600 & 0.006577965 & 0.0004359336 & -0.04101810 & 0.09227325\\
& 700 & 0.005574953 & 0.0003294627 & -0.03931742 & 0.08363312\\
& 800 & 0.005196938 & 0.0002914057 & -0.02245561 & 0.07867760\\
& 900 & 0.004681950 & 0.0002552058 & -0.03305608 & 0.07551096\\
& 1000 & 0.003797666 & 0.0002160883 & -0.01393033 & 0.06746986\\
\bottomrule
\end{tabular*}
\label{Table of Bias and MSE 1}
\end{table}

\begin{table}
\caption{Bias and MSE of MLE and moment estimator for $\mu= 0.25, 0.75, 1.5$, sample size $n=100(100)1000$, and simulation size $N=2000$}
\renewcommand{\arraystretch}{1}
\begin{tabular*}{\textwidth}{@{\extracolsep\fill}cccccc} 
\toprule
\multirow{2}{1 cm}{$\mu$} & \multirow{2}*{$n$} & \multicolumn{2}{@{}c@{}}{MLE}  & \multicolumn{2}{@{}c@{}}{Moment estimator}\\ \cmidrule{3-4}\cmidrule{5-6}%
& & Bias & MSE & Bias & MSE \\ \midrule 
\multirow{10}{1 cm}{0.25} & 100 & -0.026141601 & 0.0040926127 &  0.19908173 & 0.30299408\\
& 200 & -0.015706993 & 0.0017378035 & 0.11979148 & 0.18277572\\
& 300 & -0.010537695 & 0.0009629724 & 0.10281278 & 0.15122684\\
& 400 & -0.009388154 & 0.0006823647 & 0.07138595 & 0.11273918\\
& 500 & -0.007654403 & 0.0005241705 & 0.06670007 & 0.10613451\\
& 600 & -0.006541976 & 0.0004273222 & 0.04710248 & 0.09317236\\
& 700 & -0.005775066 & 0.0003498924 & 0.03681952 & 0.08264223\\
& 800 & -0.005424833 & 0.0003076071 & 0.02613340 & 0.07893736\\
& 900 & -0.004687041 & 0.0002641379 & 0.02308814 & 0.07288723\\
& 1000 & -0.004587965 & 0.0002188929 & 0.02168387 & 0.07098162 \\
\midrule
\multirow{10}{1 cm}{0.75} & 100 & -0.039403641 & 0.0141163154 & 0.02381151 & 0.36120747\\
& 200 & -0.023215667 & 0.0057148979 & -0.01903135 & 0.24011612\\
& 300 & -0.018010206 & 0.0037173032 & -0.04910507 & 0.18917471\\
& 400 & -0.013004298 & 0.0026953723 & -0.02435556 & 0.13988301\\
& 500 & -0.011937422 & 0.0020292707 & -0.04849263 & 0.13064004\\
& 600 & -0.010423756 & 0.0016294899 & -0.02277927 & 0.10673985\\
& 700 & -0.009229675 & 0.0013339162 & -0.03201923 & 0.09343322\\
& 800 & -0.009286969 & 0.0011661790 & -0.02381294 & 0.08053673\\
& 900 & -0.008850638 & 0.0010360152 & -0.02366172 & 0.07657251\\
& 1000 & -0.008957869 & 0.0009469547 & -0.01562601 & 0.06498910\\ 
\midrule
\multirow{10}{1 cm}{1.5} & 100 & -0.05884554 & 0.047383765 & -0.013127745 & 0.54926083\\
& 200 & -0.03000957 & 0.020339855 & -0.006090958 & 0.25400260\\
& 300 & -0.03186303 & 0.013144354 & -0.010153344 & 0.16976755\\
& 400 & -0.02284458 & 0.009625618 & -0.002482506 & 0.13057787\\
& 500 & -0.02293975 & 0.007759937 & -0.008421723 & 0.09775602\\
& 600 & -0.01841574 & 0.005969112 & -0.009321518 & 0.08234397\\
& 700 & -0.01724714 & 0.005152509 & -0.004056375 & 0.06780765\\
& 800 & -0.01197158 & 0.004282729 & -0.008411236 & 0.06338525\\
& 900 & -0.01484989 & 0.004060138 & -0.009773787 & 0.05273530\\
& 1000 & -0.01228560 & 0.003408502 & 0.003389947 & 0.04925195\\ 
\bottomrule
\end{tabular*}
\label{Table of Bias and MSE 2}
\end{table}

\section{Application} \label{sect application}
In this section, we present the application of skew-symmetric-Laplace-uniform distribution for modeling daily percentage change in the price of NIFTY 50, an Indian stock market index. Further, we have fitted and compared the proposed distribution $SSLUD(\mu)$ with normal distribution $N(\theta, \sigma^2)$, Laplace distribution $L(\theta, \beta)$, and skew-Laplace distribution $SL(\lambda)$ for percentage change data. Here, $SL(\lambda)$ refers to a special case of skew-Laplace distribution using $f$ and $K$ as pdf and cdf of standard Laplace distribution in (\ref{azzalini}). The NIFTY 50 is a benchmark Indian stock market index representing the weighted average of 50 of the largest Indian companies listed on the National Stock Exchange (NSE). It is one of the two leading stock indices used in India. The daily price of NIFTY 50 quoted in the National Stock Exchange of India Ltd. is available at \href{https://in.investing.com/indices/s-p-cnx-nifty-historical-data}{https://in.investing.com/indices/s-p-cnx-nifty-historical-data} and is selected for the current study. We consider the daily percentage change $Y_t$ on day t given by $ \displaystyle Y_t = \frac{X_t \, - \, X_{t-1}}{X_{t-1}} \times 100$, where $X_t$ denotes the price of NIFTY 50 on day t. This transformed data covering the period $16^{th}$ December 2021 to $13^{th}$ April 2022 (82 working days) is as follows :\\
0.16, -\,1.53, -\,2.18, 0.94, 1.10, 0.69, -\,0.40, 0.49, 0.86, -\,0.11, -\,0.06, 0.87, 1.57, 1.02, 0.67, -\,1.00, 0.38, 1.07, 0.29, 0.87, 0.25, -\,0.01, 0.29, -\,1.07, -\,0.96, -\,1.01, -\,0.79, -\,2.66, 0.75, -\,0.97, -\,0.05, 1.39, 1.37, 1.16, -\,1.24, -\,0.25, -\,1.73, 0.31, 1.14, 0.81, -\,1.31, -\,3.06, 3.03, -\,0.17, -\,0.10, -\,0.16, -\,0.40, -\,0.67, -\,0.17, -\,4.78, 2.53, 0.81, -\,1.12, -\,0.65, -\,1.53, -\,2.35, 0.95, 2.07, 1.53, 0.21, 1.45, -\,1.23, 1.87, 1.84, -\,0.98, 1.16, -\,0.40, -\,0.13, -\,0.40, 0.40, 0.60, 1.00, -\,0.19, 1.18, 2.17, -\,0.53, -\,0.83, -\,0.94, 0.82, -\,0.62, -\,0.82, -\,0.31.

Mean, variance, and skewness for the above data is 0.027, 1.671, and -\,0.639 respectively.
 The Wald-Wolfowitz runs test for randomness of $Y_t$ yields a p-value of 0.076, justifying the assumption of independence of the $Y_t$ values. We consider fitting the proposed skew-symmetric-Laplace-uniform distribution $SSLUD(\mu)$ along with normal distribution $N(\theta, \sigma^2)$, Laplace distribution $L(\theta, \beta)$, and skew-Laplace distribution $SL(\lambda)$ to the data on percentage change. Using R-software, the MLE of the parameters and hence, the estimated value of log-likelihood are obtained. Akaike's Information Criteria (AIC) and Bayesian Information Criteria (BIC) are used for model comparison. Table \ref{Table of AIC BIC 1} shows that the proposed $SSLUD(\mu)$ provides the best fit for the data set which is very close to $SL(\lambda)$ in terms of BIC. But in terms of AIC, N($\theta, \sigma^2$) seems to be better than $SSLUD(\mu)$ and the best among the four distributions.  

\begin{table}
\caption{MLEs, log-likelihood, AIC and BIC for daily percentage change in Nifty 50 index price ($Y_t$) of 82 days}
\begin{tabular}{@{}lcccc@{}} 
\toprule
Distribution & MLEs & ln{$L$} & AIC & BIC \\
\midrule
$SSLUD(\mu)$ & $\hat{\mu}$= 62.38674
 & -\,138.7604 & 279.5207 & 281.9274\\
 \midrule
$SL(\lambda)$ & $\hat{\lambda}$=  -6.247468e-05 & -\,138.7782 & 279.5564 & 281.9631\\
\midrule
$L(\theta, \beta)$ & $\hat{\theta}$=  -\,0.03
 , $\hat{\beta}$= 0.9990244 & -\,138.7580  
 & 281.5161 & 286.3295\\
\midrule
$N(\theta, \sigma^2)$ & $\hat{\theta}$= 0.02682927, $\hat{\sigma}^2$= 1.650275
 & -\,136.9081 & 277.8162 & 282.6296\\
\bottomrule
\end{tabular}
\label{Table of AIC BIC 1}
\end{table}

\begin{figure}
    \centering
    \includegraphics[width=0.8\textwidth]{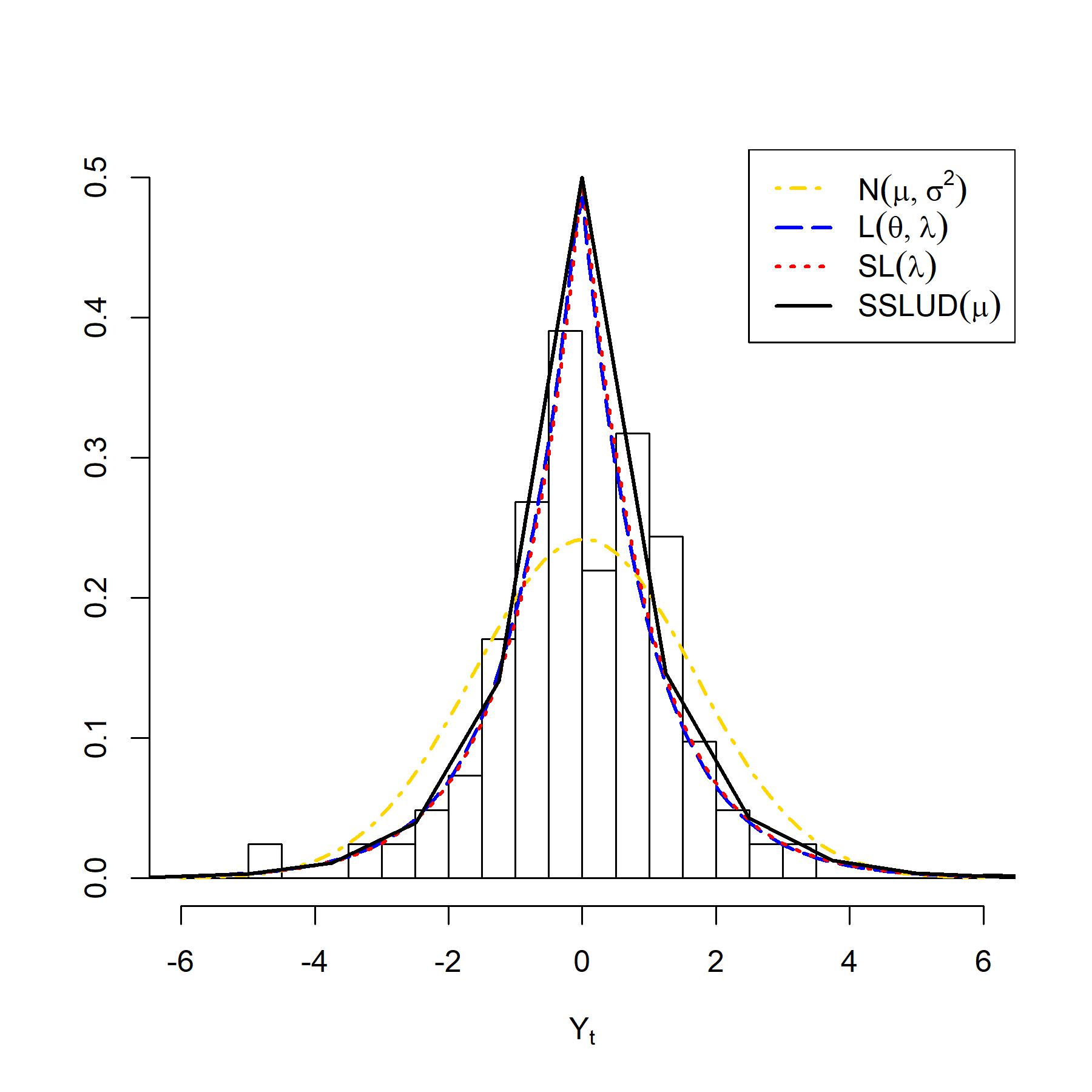}
    \caption{Plot of observed and expected densities of Normal distribution, Laplace distribution, skew-Laplace distribution, and SSLUD for daily percentage change in Nifty 50 index price ($Y_t$) of 82 days}
    \label{fig: comparison of SSLUD, Normal and Laplace distribution function 1}
\end{figure}
 
 For $SSLUD(\mu)$, MLE of $\mu$ is $\hat{\mu}$=62.38674 which is relatively high, and by definition of $g(x)$ in (\ref{pdf.SSLUD}) for a large value of $\mu$, SSLUD approaches to Laplace distribution. But, from the histogram in Figure \ref{fig: comparison of SSLUD, Normal and Laplace distribution function 1} and the value of skewness for $Y_t$, one can observe that data is negatively skewed. It might be due to a single parameter in the proposed distribution unable to give the best fit to the data. By changing the data location, significant change observed in SSLUD's curve in terms of location, scale, and shape. So, by observation, one can choose an appropriate change in location such that the value of $\hat{\mu}$ is significantly small for possible better fitting of data using the proposed distribution. Through the trial and error method, here we consider a change as - 0.8 and define transformed daily percentage change in Nifty 50 index price, $Z_t=Y_t-0.8$ which gives $\hat{\mu}= -2.589259$, a significantly small value. A possible generalization of proposed distribution with additional location parameter to avoid hindrance to employ it is under consideration, in order to make it more flexible and apt to catch the features present in real data.
 
 Table \ref{Table of AIC BIC 2} shows the MLEs, estimated log-likelihood, AIC, and BIC by fitting the distributions mentioned above to $Z_t$. The graphical representation of the results is given in Figure \ref{fig: comparison of SSLUD, Normal and Laplace distribution function 2}. It is clear from Table \ref{Table of AIC BIC 2} that the proposed $SSLUD(\mu)$ provides the best fit for the data set in terms of both AIC and BIC, but close to $SL(\lambda)$. The plot of observed and expected densities presented in Figure \ref{fig: comparison of SSLUD, Normal and Laplace distribution function 2} also confirms our findings. Thus, $SSLUD(\mu)$ is better for modeling daily percentage change in the price of NIFTY 50 in comparison to normal distribution $N(\theta, \sigma^2)$ and Laplace distribution $L(\theta, \beta)$, and one good alternative to skew-Laplace distribution $SL(\lambda)$.

\begin{table}
\caption{MLEs, log-likelihood, AIC and BIC for transformed daily percentage change in Nifty 50 index price ($Z_t$) of 82 days}
{\begin{tabular}{@{}lcccc@{}} 
\toprule
Distribution & MLEs & ln{$L$} & AIC & BIC \\
\midrule
$SSLUD(\mu)$ & $\hat{\mu}$= -\,2.589259
 & -\,136.8343 & 275.6685 & 278.0752\\
 \midrule
$SL(\lambda)$ & $\hat{\lambda}$= -\,0.6988722
 & -\,137.0020 & 276.0040 & 278.4107\\
\midrule
$L(\theta, \beta)$ & $\hat{\theta}$= -\,0.83
 , $\hat{\beta}$= 0.9990244 & -\,138.7580  
 & 281.5161 & 286.3295\\
\midrule
$N(\theta, \sigma^2)$ & $\hat{\theta}$= -\,0.7731707, $\hat{\sigma}^2$= 1.650275
 & -\,136.9081 & 277.8162 & 282.6296\\
\bottomrule
\end{tabular}}
\label{Table of AIC BIC 2}
\end{table}

\begin{figure}
    \centering
    \includegraphics[width=0.8\textwidth]{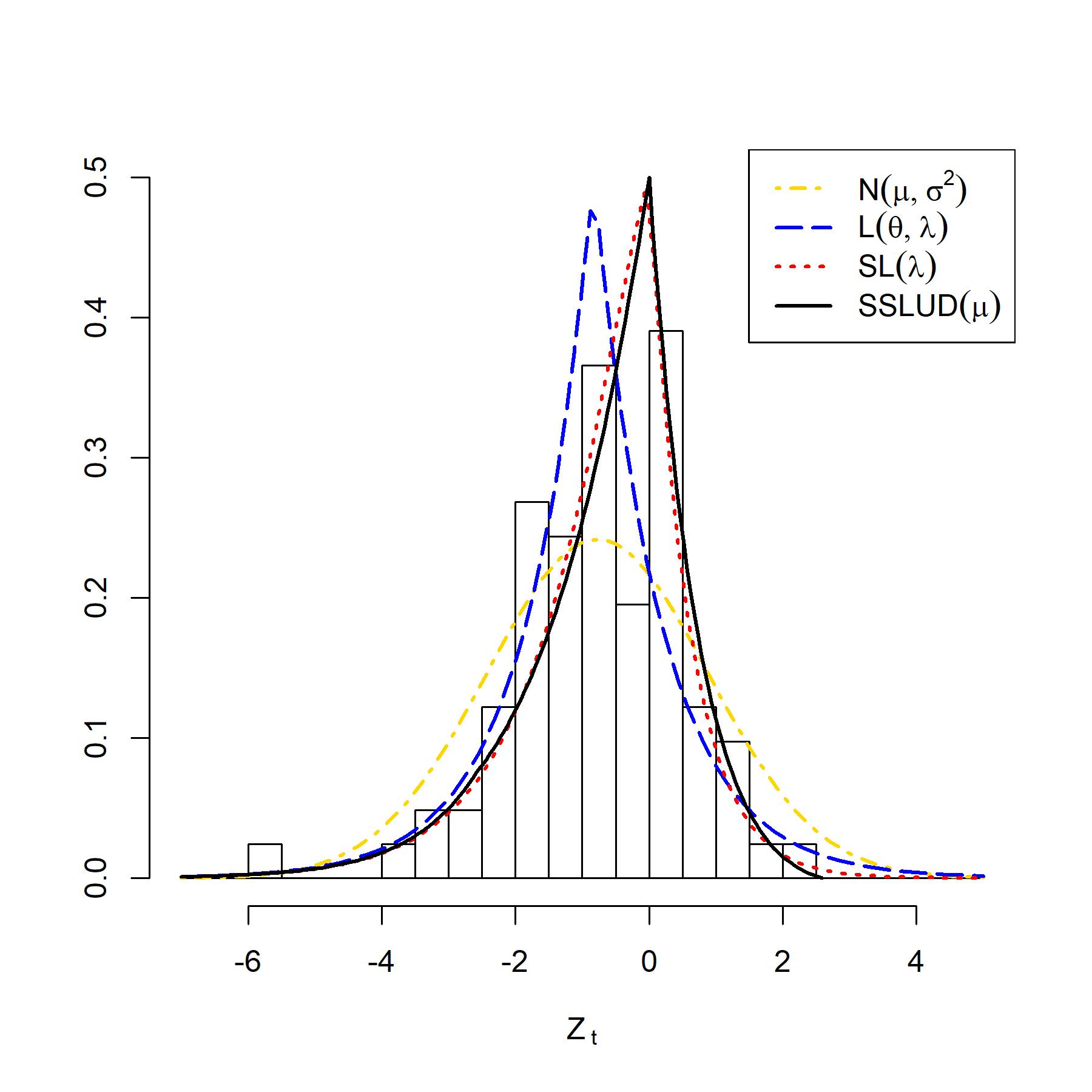}
    \caption{Plot of observed and expected densities of Normal distribution, Laplace distribution, skew-Laplace distribution, and SSLUD for transformed daily percentage change in Nifty 50 index price ($Z_t$) of 82 days}
    \label{fig: comparison of SSLUD, Normal and Laplace distribution function 2}
\end{figure}

\section*{Declarations}

\subsection*{Conflict of interest} The authors declare that they have no conflict of interest to disclose.

\subsection*{Ethics approval and consent to participate} Not applicable.

\subsection*{Consent for publication} Both authors have agreed to submit and publish this paper.

\clearpage
\bibliography{citation}

\begin{thebibliography}{}

\bibitem[Arnold and Lin, 2004]{arnold2004characterizations}
Arnold, B.~C. and Lin, G.~D. (2004).
\newblock Characterizations of the skew-normal and generalized chi distributions.
\newblock {\em Sankhy{\=a}: The Indian Journal of Statistics}, 66(4):1--14.

\bibitem[Aryal and Rao, 2005]{aryal2005reliability}
Aryal, G. and Rao, A. (2005).
\newblock Reliability model using truncated skew-laplace distribution.
\newblock {\em Nonlinear Analysis}, 63:e639--e646.

\bibitem[Azzalini, 1985]{azzalini1985class}
Azzalini, A. (1985).
\newblock A class of distributions which includes the normal ones.
\newblock {\em Scandinavian journal of statistics}, 12(2):171--178.

\bibitem[Kotz et~al., 2001]{kotz2001laplace}
Kotz, S., Kozubowski, T., and Podg{\'o}rski, K. (2001).
\newblock {\em The Laplace distribution and generalizations: a revisit with applications to communications, economics, engineering, and finance}.
\newblock Springer Science \& Business Media.

\bibitem[Kozubowski and Nolan, 2008]{kozubowski2008infinite}
Kozubowski, T.~J. and Nolan, J.~P. (2008).
\newblock Infinite divisibility of skew gaussian and laplace laws.
\newblock {\em Statistics \& probability letters}, 78(6):654--660.

\bibitem[Kozubowski and Podg{\'o}rski, 2001]{kozubowski2001asymmetric}
Kozubowski, T.~J. and Podg{\'o}rski, K. (2001).
\newblock Asymmetric laplace laws and modeling financial data.
\newblock {\em Mathematical and Computer Modelling}, 34:1003--1021.

\bibitem[Nadarajah, 2009]{nadarajah2009skew}
Nadarajah, S. (2009).
\newblock The skew logistic distribution.
\newblock {\em Advances in Statistical Analysis}, 93:187--203.

\bibitem[Nadarajah and Kotz, 2003]{nadarajah2003skewed}
Nadarajah, S. and Kotz, S. (2003).
\newblock Skewed distributions generated by the normal kernel.
\newblock {\em Statistics \& probability letters}, 65(3):269--277.

\bibitem[Nekoukhou and Alamatsaz, 2012]{nekoukhou2012family}
Nekoukhou, V. and Alamatsaz, M. (2012).
\newblock A family of skew-symmetric-laplace distributions.
\newblock {\em Statistical papers}, 53:685--696.

\bibitem[Rachev and SenGupta, 1993]{rachev1993laplace}
Rachev, S. and SenGupta, A. (1993).
\newblock Laplace-weibull mixtures for modeling price changes.
\newblock {\em Management Science}, 39(8):1029--1038.

\bibitem[Ryd{\'e}n et~al., 1998]{ryden1998stylized}
Ryd{\'e}n, T., Ter{\"a}svirta, T., and {\AA}sbrink, S. (1998).
\newblock Stylized facts of daily return series and the hidden markov model.
\newblock {\em Journal of applied econometrics}, 13(3):217--244.

\bibitem[Theodossiou, 1998]{theodossiou1998financial}
Theodossiou, P. (1998).
\newblock Financial data and the skewed generalized t distribution.
\newblock {\em Management Science}, 44(Part 1 of 2):1650--1661.

\bibitem[Zeckhauser and Thompson, 1970]{zeckhauser1970linear}
Zeckhauser, R. and Thompson, M. (1970).
\newblock Linear regression with non-normal error terms.
\newblock {\em The Review of Economics and Statistics}, 52:280--286.

\end{thebibliography}
\bibliographystyle{apalike}

\end{document}